\documentclass[11pt, reqno]{amsart}
\usepackage{amsthm,amsmath,amsfonts,amssymb,color}
\usepackage[bookmarks]{hyperref}
\reversemarginpar

\newtheorem{thm}{Theorem}[section]

\newtheorem{preremark}[thm]{Remark}
\newenvironment{remark}{\begin{preremark}\rm}{\medskip \end{preremark}}

\numberwithin{equation}{section}

\newcommand{\abs}[1]{\left\vert#1\right\vert}

\newcommand{\R}{\mathbb R}

\newcommand{\dd} {\mathrm{d}}

\DeclareMathOperator{\dv}{div}
\DeclareMathOperator{\Def}{Def}
\DeclareMathOperator{\Ric}{Ric}

\def\be{\begin{equation}}
\def\ee{\end{equation}}
\def\bs{\begin{split}}
\def\ess{\end{split}}
\def\o{\omega}

\begin{document}
\title[Ellipsoid]{The restriction problem on the ellipsoid}

\author[Chan]{Chi Hin Chan}
\address{Department of Applied Mathematics, National Yang Ming Chiao Tung University,1001 Ta Hsueh Road, Hsinchu, Taiwan 30010, ROC}
\email{cchan@math.nctu.edu.tw}

\author[Czubak]{Magdalena Czubak}
\address{Department of Mathematics\\
University of Colorado Boulder\\ Campus Box 395, Boulder, CO, 80309, USA}
\email{czubak@math.colorado.edu}
%
%
\textcolor{red}{
\author[Yoneda]{Tsuyoshi Yoneda} 
\address{Graduate School of Economics, Hitotsubashi University, 2-1 Naka, Kunitachi, Tokyo 186-8601, Japan} 
\email{t.yoneda@r.hit-u.ac.jp} 
}

\begin{abstract}
Following a restriction argument in the Euclidean space, we derive a geometric invariant formula for a possible viscosity operator for an incompressible fluid flow on an ellipsoid embedded in $\R^3$. We also give an asymptotic expansion of the formula in terms of the eccentricity associated with the ellipsoid.
\end{abstract}
\maketitle

 \section{Introduction}
This is a continuation of work started in \cite{CCM17}.  There the first two authors and Marcelo Disconzi investigated the question of what should the viscosity operator be on a Riemannian manifold for divergence free vector fields.  In the Euclidean case, the viscosity operator is the Laplacian.  However, there are several different choices of the Laplacian that can be applied to a vector field on a manifold.  These include the Bochner Laplacian and the Hodge Laplacian.

In their 1970 article \cite{EbinMarsden}, Ebin and Marsden indicated that
when
writing the Navier-Stokes equation on an Einstein
manifold, one should use the following operator
\be\label{deformation}
 2\Def^\ast \Def,
\ee
where $\Def$ is the deformation tensor, and $\Def^\ast$ is its adjoint.  In coordinates, we can write $\Def$ as
\begin{gather}
 {(\Def u)}_{ij}=\frac 12 (\nabla_i u_j +\nabla_j u_i),
\end{gather}
where $\nabla$ is the Levi-Civita connection. In the Euclidean case, operator in \eqref{deformation} reduces to the (standard) Laplacian.  Similarly,  the Bochner Laplacian and the Hodge Laplacian coincide, and are also equal to the (standard) Laplacian.  On a general manifold, these operators are not the same.

 In \cite{CCM17}, there was further evidence provided why one might want to follow Ebin and Marsden.  One argument came from the so-called restriction argument: start with a divergence free vector field on $\R^3$, apply the Euclidean Laplacian to the vector field, and then restrict to a sphere, and see what operator one gets in that case.  In the case of the sphere, the computation in \cite{CCM17} shows that it is indeed the operator in \eqref{deformation}.
 
We would like to note that there are other works which address the question of the formulation of the Navier-Stokes equations on manifolds.  See for example \cite{Kobayashi, Fang2020, Samavaki_Tuomela2020} and the references therein.

In this paper, we continue this line of investigation by looking at the restriction argument as applied to an ellipsoid embedded in $\R^3$.  One motivation is to consider a space that is not a space form, but still has some symmetry.  Another motivation is that the ellipsoid can be viewed as a model for Jupiter.

The incompressible 2D-Navier-Stokes equation on a rotating sphere 
is one of the simplest models of planets such as Jupiter and Earth, and many researchers have been extensively studying this model. Williams \cite{W} was the first to find that 
turbulent flow becomes a multiple jet flow on such a model. 
After that Yoden-Yamada \cite{YY}, Nozawa-Yoden \cite{NY}, Obuse-Takehiro-Yamada \cite{O} and Sasaki-Takehiro-Yamada \cite{STY-1,STY-2} made further progress.
%

For the viscosity operator, $-\Delta-2$ is the standard formula for that model. The additional term $``2"$ (this is corresponding to the  Ricci curvature)
 is necessary for the conservation of the total angular momentum of the system \cite{STY-2}.  
 
However, as far as the authors are aware, 
none of the numerous works on multiple jet flow model attempted to investigate the effect of the ``underlying manifold'' itself on the viscosity effect.  
More precisely,  in the above simplest model, the manifold is a ``sphere", even though in reality
Jupiter and Earth are not spheres. They have ``bulges" around the equator and are flattened at the poles, and therefore, they could be viewed as ellipsoids. 
 
To derive the appropriate viscosity operator on the ellipsoid,
 we need to clarify what kind of 
approximation we need to apply (c.f. in the inviscid flow case, we do not need such approximation. See for example \cite{TY1, TY2} for the study of the inviscid flow on the ellipsoid).
%

Physically, centrifugal force affects the atmosphere proportionally to the distance from the axis of rotation.  Thus this simple physical observation naturally supports the use of the following approximation of the atmospheric layer: 
\begin{equation}\label{parametrized}
\Phi(\rho, \phi, \theta)=(a\rho\sin \phi \cos\theta, a\rho \sin\phi \sin\theta, \rho\cos \phi)
\end{equation}
with the parameter $\rho\geq 1$ (c.f. \cite{Miura}).

There are some difficulties in carrying out the computations in the case of the ellipsoid.  First, the metric is more complicated.  Second,  there are difficulties related to the divergence free condition.  Vector fields in $\R^3$ that restrict to vector fields on a sphere and are divergence free and independent of the radius in $\R^3$ are also divergence free on the sphere. This is no longer the case on the ellipsoid.  To address this, we work with the most general vector field on $\R^3$ that is also tangential to the ellipsoid, and we assume that the vector field satisfies both of the divergence free conditions: one on $\R^3$, and one on the ellipsoid. Such vector fields can be shown to exist (see Section \ref{divfree}).

 Consider the ellipsoid defined by
\[
E=\{ x^2+y^2+a^2z^2=a^2\},
\]
where $a>0$.  The ellipsoid $E$ can be parametrized by a map
\[
\Phi:  (0,\infty)\times (0,\pi)\times (-\pi,\pi) \to \R^3
\]
given by \eqref{parametrized}
when $\rho=1$.

We are now ready to state the main result.

 \begin{thm}\label{main}
 Let $v^\sharp$ be a vector field defined in a neighborhood of the ellipsoid $E$, that is tangential to the ellipsoid, and that is both divergence free on $\R^3$ and on the ellipsoid\footnote{The precise formula for $v$ is given in the section below.}.  Then

\begin{equation}\label{formula}
\begin{split}
\mathfrak{i}_{\mathrm{E}}^* \Big \{ -\triangle v \Big \} =   -\triangle_{\mathrm{E}}\Big ( \mathfrak{i}_{\mathrm{E}}^* v  \Big ) + \mathcal{E} \big ( v \big )
- \sqrt{K_{\mathrm{E}}}  \mathfrak{i}_{\mathrm{E}}^* \Big \{ \mathcal{L}_{\mathrm{Y}} \big ( v \big ) \Big \}
- 2 \sqrt[4]{K_{\mathrm{E}}} \Big (1- \frac{1}{|\nabla \rho |^2} \Big )\mathfrak{i}_{\mathrm{E}}^* \Big \{ (\mathcal L_{\nabla \rho} v)_\phi \} \Big \} e^2,\\
\end{split}
\end{equation}
where
\begin{itemize}
\item $v=(v^\sharp)^\flat$, where $\sharp$ and $\flat$ are the standard musical isomorphisms;
\item $\mathfrak i^\ast_E$ is the pullback by the inclusion map $\mathfrak i_E: E\hookrightarrow \R^3$;
\item $-\triangle$ is the Hodge Laplacian on $\R^3$;
\item  $ -\triangle_{\mathrm{E}}$ is the Hodge Laplacian on $E$;
\item $\mathcal{E}$ is an operator given by
\begin{equation}\label{operator}
\mathcal{E} = \mathfrak{i}_{\mathrm{E}}^* \Big \{ -\mathcal{L}_{\nabla \rho} \Big ( \frac{1}{|\nabla \rho|^2}  \mathcal{L}_{\nabla \rho}     \Big ) + \Big ( 1 - \frac{1}{|\nabla \rho|^2}     \Big )\mathcal{L}_{\nabla \rho} \Big \} + \sqrt{\mathrm{K}_{\mathrm{E}}}  \mathfrak{i}_{\mathrm{E}}^* \Big \{ \frac{1}{|\nabla \rho|^2}  \mathcal{L}_{\nabla \rho}  \Big \} ,
\end{equation}
where $\mathcal L$ denotes the Lie derivative, $K_E$ is the sectional curvature of the ellipsoid;
\item $Y=\rho \partial_\rho$;
\item $e^2$ is the $1$-form on $E$ dual to  $e_2$, which is a unit vector field in the direction of $\partial_\phi$, the coordinate vector field that points along the latitude directions on the ellipsoid;

\end{itemize}

\end{thm}

\begin{remark}
For divergence free vector fields/forms, the operator in \eqref{deformation} is related to the Hodge Laplacian,  $\dd\dd^\ast +\dd^\ast \dd$, by

 \be\label{divdeffull}
 2\Def^\ast \Def=\dd\dd^\ast +\dd^\ast \dd -2\Ric,
\ee
where $\Ric$ is the Ricci curvature tensor (with a raised index).

In \cite{CCM17}, the restriction argument produced exactly the Hodge Laplacian, $\dd\dd^\ast +\dd^\ast \dd$, and  $-2\Ric$.  In the case of the ellipsoid, we obtain the Hodge Laplacian and several other terms, without the Ricci term explicitly appearing.
\end{remark}
\begin{remark}
One can show that in the case of the sphere the above formula does reduce to what was obtained in \cite{CCM17}.  See Section \ref{sphere}.
\end{remark}

In the next section, we prove Theorem \ref{main}.  In Section \ref{sphere}, as indicated above, we show how the formula simplifies in the case of the sphere.  Section \ref{divfree} is dedicated to showing how one can construct a vector field that is simultaneously divergence free on the ellipsoid and on $\R^3$.  Finally, in Section \ref{power_exp}, we give an asymptotic expansion of the formula in terms of the eccentricity associated with the ellipsoid.

\section{Set-up and Proof of Theorem \ref{main}}\label{proof}

Let
\[
E=\{ x^2+y^2+a^2z^2=a^2\},
\]
where $a>0$, and let $E$ be parametrized by the map
\[
\Phi:  (0,\infty)\times (0,\pi)\times (-\pi,\pi) \to \R^3
\]
given by
\[
\Phi(\rho, \phi, \theta)=(a\rho\sin \phi \cos\theta, a\rho \sin\phi \sin\theta, \rho\cos \phi),
\]
when $\rho=1$.  The metric on the ellipsoid is the Euclidean metric restricted to $E$.  We begin with the metric on $\R^3$. In these coordinates, the Euclidean metric on $\R^3$ is given by
\be\label{Rmetric}
\left(\begin{array}{ccc}
   g_{\rho \rho}&g_{\rho\phi}&0\\
   g_{\phi\rho }&g_{\phi\phi}&0\\
 0&0&g_{\theta\theta}\\
   \end{array}\right),
\ee
where
\begin{align}
g_{\rho\rho}&=a^2 \sin^2\phi+\cos^2\phi, \quad\qquad g_{\rho \phi}=g_{\phi\rho}=(a^2-1)\rho \sin\phi\cos\phi,\\
g_{\phi \phi}&=a^2\rho^2\cos^2\phi+\rho^2\sin^2\phi, \quad g_{\theta\theta}=a^2\rho^2\sin^2\phi.
\end{align}
  
The metric on $E$ is
\be\label{Emetric}
\left(\begin{array}{cc}
   g_{\phi \phi}&0\\
   0&g_{\theta\theta}\\
   \end{array}\right),
\ee
with $g_{\phi\phi}, g_{\theta\theta}$ as above, with $\rho=1$.

Let $\star$ denote the Hodge star operator.  We use the following
\begin{align}
\star \dd\rho \wedge \dd\phi&=\sin \phi \dd\theta,\\
\star \dd \rho\wedge \dd \theta&=\frac{(1-a^2)\cos \phi}{a^2\rho} \dd \rho-\frac{\lambda^2}{a^2\sin\phi} \dd\phi,\\
\star \dd\phi\wedge \dd \theta &= \frac{\tfrac{(1-a^2)^2}{a^2}\sin^2\phi\cos^2\phi+1}{\rho^2 \lambda^2\sin\phi}\dd\rho-\frac{(1-a^2)\cos\phi}{a^2\rho}\dd \phi,
\end{align}
where $\lambda^2=a^2\cos^2\phi+\sin^2\phi$.
Let $v^{\sharp}$ be a smooth vector field on $\R^3$.  In above coordinates
\[
v^{\sharp}=v^{\rho}\partial_\rho+v^\phi\partial_\phi+v^\theta\partial_\theta.
\]
We also let $v^\rho=f(\rho, \phi, \theta)(\rho-1)$ so that $v^\rho=0$ on $E$, and suppose
\[
\dv_{\R^3}v^\sharp=0=\dv_{E}v^\sharp.
\]
This means
\[
\dv_{\R^3}v^\sharp=\partial_\rho v^{\rho}+\partial_\phi v^\phi +\frac 2\rho v^\rho+\partial_\theta v^\theta +\cot\phi v^\phi=0,
\]
and
\[
\dv_{E}v^\sharp=\partial_\phi v^\phi +\partial_\theta v^\theta +\cot\phi\frac{(2-a^2)\sin^2\phi+a^2\cos^2\phi}{a^2\cos^2\phi+\sin^2\phi} v^\phi=0.
\]

Next, the $1$-form $v$ is
\begin{align*}
v&=(g_{\rho \rho}v^\rho+g_{\rho \phi}v^\phi) \dd\rho+g_{\phi \phi}v^\phi \dd\phi+g_{\theta\theta}v^\theta \dd \theta\\
&=: v_\rho \dd\rho+v_\phi \dd\phi+v_\theta \dd \theta,
\end{align*}
so that
\[
\omega=\dd v=\omega_{\rho\phi}\dd\rho\wedge\dd\phi +\omega_{\rho\theta}\dd\rho\wedge\dd\theta+\omega_{\phi\theta}\dd\phi\wedge\dd\theta,
\]
where
\[
\omega_{ij}=\partial_i v_j-\partial_j v_i.
\]
It follows
\begin{align*}
\star\dd v&=\omega_{\rho\phi} \sin\phi \dd\theta +\omega_{\rho\theta}(\frac{(1-a^2)\cos \phi}{a^2\rho} \dd \rho-\frac{\lambda^2}{a^2\sin\phi} \dd\phi)\\
&\qquad+\omega_{\phi\theta}(\frac{\tfrac{(1-a^2)^2}{a^2}\sin^2\phi\cos^2\phi+1}{\rho^2 \lambda^2\sin\phi}\dd\rho-\frac{(1-a^2)\cos\phi}{a^2\rho}\dd \phi)\\
&=\left(\omega_{\rho\theta} \frac{(1-a^2)\cos \phi}{a^2\rho}  +\omega_{\phi\theta}\frac{\tfrac{(1-a^2)^2}{a^2}\sin^2\phi\cos^2\phi+1}{\rho^2 \lambda^2\sin\phi} \right)\dd\rho-\left( \omega_{\rho\theta}\frac{\lambda^2}{a^2\sin\phi} +\omega_{\phi\theta}\frac{(1-a^2)\cos\phi}{a^2\rho} \right)\dd\phi\\
&\qquad
+\omega_{\rho\phi} \sin\phi \dd\theta\\
&=:A_\rho \dd\rho +A_\phi \dd\phi +A_\theta \dd\theta.
\end{align*}
Then, similarly as above
\[
\dd\star \dd v=F_{\rho\phi}\dd\rho\wedge\dd\phi +F_{\rho\theta}\dd\rho\wedge\dd\theta+F_{\phi\theta}\dd\phi\wedge\dd\theta,
\]
with
\[
F_{ij}=\partial_i A_j-\partial_j A_i.
\]
Using that
\[
-\Delta v=\dd^\star \dd v= \dd\star \dd v,
\]
we obtain
\be
\begin{split}
-\Delta v&=\left(F_{\rho\theta} \frac{(1-a^2)\cos \phi}{a^2\rho}  +F_{\phi\theta}\frac{\tfrac{(1-a^2)^2}{a^2}\sin^2\phi\cos^2\phi+1}{\rho^2 \lambda^2\sin\phi} \right)\dd\rho-\left( F_{\rho\theta}\frac{\lambda^2}{a^2\sin\phi} +F_{\phi\theta}\frac{(1-a^2)\cos\phi}{a^2\rho} \right)\dd\phi\\
&\qquad+F_{\rho\phi} \sin\phi \dd\theta
\end{split}
\ee
Next
\be\label{lhs}
\mathfrak{i}_{\mathrm{E}}^* \Big \{ -\Delta v \Big \} =\mathfrak{i}_{\mathrm{E}}^*\Big\{ - \left( F_{\rho\theta}\frac{\lambda^2}{a^2\sin\phi} +F_{\phi\theta}\frac{(1-a^2)\cos\phi}{a^2} \right)\dd\phi
+F_{\rho\phi} \sin\phi \dd\theta \Big\}.
\ee
To show \eqref{formula} holds, we consider the right hand side and compute all the terms and then compare with \eqref{lhs}.
\subsection{Hodge Laplacian on $E$}
Since $v^\sharp$ is divergence free, we just need to consider $\dd^\star \dd$.  A computation shows
\be
-\Delta_E v=\frac{\partial_\theta \o_{\phi\theta}}{a^2\sin^2\phi}\dd\phi-\frac{\sin \phi}{\lambda}\partial_\phi \left(\frac{\o_{\phi\theta}}{\lambda\sin\phi}\right)\dd\theta.
\ee
\subsection{The operator $\mathcal E(v)$}
We recall

\begin{equation}
\mathcal{E} = \mathfrak{i}_{\mathrm{E}}^* \Big \{ -\mathcal{L}_{\nabla \rho} \Big ( \frac{1}{|\nabla \rho|^2}  \mathcal{L}_{\nabla \rho}     \Big ) + \Big ( 1 - \frac{1}{|\nabla \rho|^2}     \Big )\mathcal{L}_{\nabla \rho} \Big \} + \sqrt{\mathrm{K}_{\mathrm{E}}}  \mathfrak{i}_{\mathrm{E}}^* \Big \{ \frac{1}{|\nabla \rho|^2}  \mathcal{L}_{\nabla \rho}  \Big \}.
\end{equation}
Here we use the Cartan formula formula for the Lie derivative of a $1$-form $\alpha$, which reads
\be
\mathcal L_X\alpha=\iota_X\dd \alpha+ \dd(\iota_X \alpha),
\ee
and $\iota$ is an interior multiplication.  Note
\[
\nabla \rho=\frac{\lambda^2}{a^2}\partial_\rho+\frac{1-a^2}{a^2\rho}\sin\phi\cos\phi\partial_\phi,
\]
and
\begin{align}
 \iota_{\nabla \rho}(\dd\rho\wedge \dd\phi)&=\frac{\lambda^2}{a^2}\dd \phi+\frac{a^2-1}{a^2\rho}\sin\phi\cos\phi\dd\rho,\\
 \iota_{\nabla \rho}(\dd\rho\wedge \dd\theta)&=\frac{\lambda^2}{a^2}\dd \theta,\\
 \iota_{\nabla \rho}(\dd\phi\wedge \dd\theta)&=\frac{1-a^2}{a^2\rho}\sin\phi\cos\phi\dd \theta,
\end{align}
and
\[
\iota_{\nabla \rho} v=v^\rho,
\]
so
\be\label{Lie1}
\mathcal L_{\nabla \rho}v=\o_{\rho\phi}\frac{a^2-1}{a^2\rho}\sin\phi\cos\phi\dd\rho+\o_{\rho\phi} \frac{\lambda^2}{a^2}\dd \phi +(\o_{\rho\theta} \frac{\lambda^2}{a^2}+\o_{\phi\theta}\frac{1-a^2}{a^2\rho}\sin\phi\cos\phi)\dd \theta+ \dd v^\rho,
\ee
and
\begin{align}
 \mathfrak{i}_{\mathrm{E}}^* \Big \{   -\mathcal{L}_{\nabla \rho} \Big ( \frac{1}{|\nabla \rho|^2}  \mathcal{L}_{\nabla \rho}  v\Big)   \Big\}
 =\mathfrak{i}_{\mathrm{E}}^* \Big \{ G_\phi \dd\phi + G_\theta \dd\theta -\dd(  \partial_\rho v^\rho)\Big\}\nonumber,
\end{align}
where we let
\begin{align}
G_\phi&=\left(-\partial_\rho( \o_{\rho\phi}+\frac{a^2}{\lambda^2}\partial_\phi v^\rho)+\partial_\phi(\frac{a^2}{\lambda^2}\partial_\rho v^\rho-\frac{1-a^2}{\lambda^2\rho}\sin\phi\cos\phi \o_{\rho\phi})\right)\frac {\lambda^2}{a^2},
\\
G_\theta&=\left(-\partial_\rho(\o_{\rho\theta}+\o_{\phi\theta}\frac{1-a^2}{\lambda^2\rho}\sin\phi\cos\phi+\frac {a^2}{\lambda^2}\partial_\theta v^\rho)+\partial_\theta(\frac{a^2}{\lambda^2}\partial_\rho v^\rho-\frac{1-a^2}{\lambda^2\rho}\sin\phi\cos\phi \o_{\rho\phi})\right)\frac {\lambda^2}{a^2}\nonumber\\
&\quad+
\left(-\partial_\phi(\o_{\rho\theta}+\o_{\phi\theta} \frac{1-a^2}{\lambda^2\rho}\sin\phi\cos\phi)+\partial_\theta \o_{\rho\phi}\right)\frac{1-a^2}{a^2\rho}\sin\phi\cos\phi.\label{G2}
\end{align}

Using product rule, we make some simplifications as follows
\begin{align}
G_\phi=-\partial_\rho \o_{\rho\phi}\frac{\lambda^2}{a^2} -\frac{a^2}{\lambda^2}\partial_\rho v^\rho \partial_\phi(\frac {\lambda^2}{a^2}) -\partial_\phi(\frac{1-a^2}{a^2\rho}\sin\phi\cos\phi \o_{\rho\phi})+\frac{1-a^2}{\lambda^2\rho}\sin\phi\cos\phi \o_{\rho\phi}\partial_\phi(\frac{\lambda^2}{a^2}).\label{G1}
\end{align}

For convenience, we also write down
\[
 \Big ( 1 - \frac{1}{|\nabla \rho|^2}     \Big ) \mathfrak{i}_{\mathrm{E}}^* \Big \{ \mathcal{L}_{\nabla \rho}  v\Big\}=\frac{\lambda^2-a^2}{\lambda^2}\left(\o_{\rho\phi} \frac{\lambda^2}{a^2}\dd \phi +(\o_{\rho\theta} \frac{\lambda^2}{a^2}+\o_{\phi\theta}\frac{1-a^2}{a^2\rho}\sin\phi\cos\phi)\dd \theta\right),
\]
and
\[
 \sqrt{\mathrm{K}_{\mathrm{E}}}  \mathfrak{i}_{\mathrm{E}}^* \Big \{ \frac{1}{|\nabla \rho|^2}  \mathcal{L}_{\nabla \rho}  v\Big\}=\frac{a^2}{\lambda^4}\left(\o_{\rho\phi} \frac{\lambda^2}{a^2}\dd \phi +(\o_{\rho\theta} \frac{\lambda^2}{a^2}+\o_{\phi\theta}\frac{1-a^2}{a^2\rho}\sin\phi\cos\phi)\dd \theta \right).
\]
For future reference, we now breakdown   $\mathcal E$ into $\dd \phi$ and $\dd \theta$ components. We have
\begin{align}
\mathcal E(v)_\phi&=G_\phi-\partial_\phi \partial_\rho v^\rho +(\frac{\lambda^2-a^2}{a^2}+\frac{1}{\lambda^2})\o_{\rho\phi},\\
\mathcal E(v)_\theta&=G_\theta-\partial_\theta \partial_\rho v^\rho+(\frac{\lambda^2-a^2}{\lambda^2}+\frac{a^2}{\lambda^4})( \o_{\rho\theta} \frac{\lambda^2}{a^2}+\o_{\phi\theta}\frac{1-a^2}{a^2\rho}\sin\phi\cos\phi).
\end{align}

\subsection{ The  $\mathcal L_Y(v)$ term}
Using the Cartan formula as above we obtain
\be
\mathcal L_Y v= \rho \o_{\rho\phi} \dd\phi +\rho \o_{\rho\theta}\dd\theta+\dd(\rho v_\rho),
\ee
so 
\[
- \sqrt{K_{\mathrm{E}}}  \mathfrak{i}_{\mathrm{E}}^* \Big \{ \mathcal{L}_{\mathrm{Y}} \big ( v \big ) \Big \}=-\frac{1}{\lambda^2}( \o_{\rho\phi} \dd\phi +\o_{\rho\theta}\dd\theta+\dd_E( v_\rho)).
\]
\subsection{The last term}
Here the term is
\be\label{last}
- 2 \sqrt[4]{K_{\mathrm{E}}} \Big (1- \frac{1}{|\nabla \rho |^2} \Big )\mathfrak{i}_{\mathrm{E}}^* \Big \{ (\mathcal L_{\nabla \rho} v)_\phi \} \Big \} e^2=-\frac{2}{a^2}(1-\frac{a^2}{\lambda^2})\o_{\rho\phi} \dd \phi
\ee
 
We are now ready to compare the coefficient functions on both sides of the formula.  Of course, the reader can choose to do this by themselves, but since the calculation can take some time, below, we set up what we hope can be helpful.
\subsection{Comparing the $\dd\phi$ components}
On the left hand side, from \eqref{lhs}, we have
 \[
 -\left( F_{\rho\theta}\frac{\lambda^2}{a^2\sin\phi} +F_{\phi\theta}{\frac{1-a^2}{a^2}\cos\phi}\right)=I+II,
 \]
 with
 \begin{align*} 
 I&=  -\left(\partial_\rho(\o_{\rho\phi}\sin \phi)- \partial_\theta\left(\omega_{\rho\theta}{\frac{1-a^2}{a^2}\cos\phi} +\omega_{\phi\theta}\frac{\tfrac{(1-a^2)^2}{a^2}\sin^2\phi\cos^2\phi+1}{\lambda^2\sin \phi}\right)\right)\frac{\lambda^2}{a^2\sin\phi}, \\
II&=-\Big(\partial_\phi  (\o_{\rho\phi}\sin \phi)+\partial_\theta( \omega_{\rho\theta}\frac{\lambda^2}{a^2\sin\phi} +\omega_{\phi\theta}{\frac{1-a^2}{a^2}\cos\phi})\Big){\frac{1-a^2}{a^2}\cos\phi}.
 \end{align*}
 We now observe that we have a cancellation between the second terms in $I$ and $II$ as well as the third terms, so we can redefine $I$ and $II$ to be
  \begin{align}
 I&=  -\left(\partial_\rho(\o_{\rho\phi}\sin \phi)- \partial_\theta\left(\omega_{\phi\theta}\frac{1}{\lambda^2\sin \phi}\right)\right)\frac{\lambda^2}{a^2\sin\phi}, \\
II&=-\partial_\phi  (\o_{\rho\phi}\sin \phi){\frac{1-a^2}{a^2}\cos\phi}.
 \end{align}
Next, we see that the second term in $I$ is exactly equal to the $\dd\phi$ component of the Hodge Laplacian on $E$.  This leaves us with only two terms on the left hand side
\be\label{left1}
-\partial_\rho \o_{\rho\phi}\frac{\lambda^2}{a^2} -\partial_\phi  (\o_{\rho\phi}\sin \phi){\frac{1-a^2}{a^2}\cos\phi}.
\ee
We now gather the $\dd \phi$ components from the right hand side.  The Hodge Laplacian was just handled. This leaves
\begin{align*}
&G_\phi -\partial_\phi \partial_\rho v^\rho+(\frac{\lambda^2-a^2}{a^2}+\frac{1}{\lambda^2})\o_{\rho\phi} -\frac{1}{\lambda^2}\left( \o_{\rho\phi}+\partial_\phi v_\rho \right)
-2\frac{\lambda^2-a^2}{a^2\lambda^2}\o_{\rho\phi} \\
&\quad=   G_\phi  -\partial_\phi \partial_\rho v^\rho+(\frac{\lambda^2-a^2-2}{a^2}+\frac 2{\lambda^2})\o_{\rho\phi}-\frac{1}{\lambda^2}\partial_\phi v_\rho\\
&\quad =-\partial_\rho \o_{\rho\phi}\frac{\lambda^2}{a^2} -\frac{a^2}{\lambda^2}\partial_\rho v^\rho \partial_\phi(\frac {\lambda^2}{a^2}) -\partial_\phi(\frac{1-a^2}{a^2\rho}\sin\phi\cos\phi \o_{\rho\phi})+\frac{1-a^2}{\lambda^2\rho}\sin\phi\cos\phi \o_{\rho\phi}\partial_\phi(\frac{\lambda^2}{a^2})\\
&\qquad-\partial_\phi\partial_\rho v^\rho+(\frac{\lambda^2-a^2-2}{a^2}+\frac 2{\lambda^2})\o_{\rho\phi}-\frac{1}{\lambda^2}\partial_\phi v_\rho
\end{align*}
We can observe that the first term in \eqref{left1} matches with the first term in the last equation.  This means we need to show
\be\label{blah}
 -\partial_\phi  (\o_{\rho\phi}\sin \phi){\frac{1-a^2}{a^2\rho}\cos\phi}=-\partial_\phi \o_{\rho\phi}\frac{1-a^2}{a^2\rho}\sin\phi\cos\phi -\o_{\rho\phi}\frac{(1-a^2)\cos^2\phi}{a^2}
\ee

is equal to
\begin{align*}
&-\frac{a^2}{\lambda^2}\partial_\rho v^\rho \partial_\phi(\frac {\lambda^2}{a^2}) -\partial_\phi(\frac{1-a^2}{a^2\rho}\sin\phi\cos\phi \o_{\rho\phi})+\frac{1-a^2}{\lambda^2\rho}\sin\phi\cos\phi \o_{\rho\phi}\partial_\phi(\frac{\lambda^2}{a^2})\\
&\qquad -\partial_\phi\partial_\rho v^\rho+ (\frac{\lambda^2-a^2-2}{a^2}+\frac 2{\lambda^2})\o_{\rho\phi}-\frac{1}{\lambda^2}\partial_\phi v_\rho
\end{align*}
  Using product rule in the second term in the above expression, we immediately see the first term on the right hand side of  \eqref{blah}.

  This leaves us with having to show
  \begin{align*}
-  \omega_{\rho\phi}\frac{(1-a^2)\cos^2\phi}{a^2}&=-\frac{a^2}{\lambda^2}\partial_\rho v^\rho \partial_\phi(\frac {\lambda^2}{a^2}) -\partial_\phi(\frac{1-a^2}{a^2\rho}\sin\phi\cos\phi )\o_{\rho\phi}+\frac{1-a^2}{\lambda^2\rho}\sin\phi\cos\phi \o_{\rho\phi}\partial_\phi(\frac{\lambda^2}{a^2})\\
&\qquad-\partial_\phi\partial_\rho v^\rho+ (\frac{\lambda^2-a^2-2}{a^2}+\frac 2{\lambda^2})\o_{\rho\phi}-\frac{1}{\lambda^2}\partial_\phi v_\rho,
\end{align*}
which one can check by direct computation.

The comparison of the  $\dd\theta$ components is similar, and in some ways easier, and we leave it to the reader.

 \section{Reduction to the sphere}\label{sphere}
In this section we show that the formula \eqref{formula} in Theorem \ref{main}, if $E=S^2$, reduces to 
\be
\mathfrak{i}_{\mathrm{S^2}}^* \Big \{ -\triangle v \Big \} =   -\triangle_{\mathrm{S^2}} v   -  2a^2 v,
 \ee
which is what was obtained in \cite{CCM17}.

For convenience, we now recall the formula
\begin{equation}\label{formula_again}
\begin{split}
\mathfrak{i}_{\mathrm{E}}^* \Big \{ -\triangle v \Big \} =   -\triangle_{\mathrm{E}}\Big ( \mathfrak{i}_{\mathrm{E}}^* v  \Big ) + \mathcal{E} \big ( v \big )
- \sqrt{K_{\mathrm{E}}}  \mathfrak{i}_{\mathrm{E}}^* \Big \{ \mathcal{L}_{\mathrm{Y}} \big ( v \big ) \Big \}
- 2 \sqrt[4]{K_{\mathrm{E}}} \Big (1- \frac{1}{|\nabla \rho |^2} \Big )\mathfrak{i}_{\mathrm{E}}^* \Big \{ (\mathcal L_{\nabla \rho} v)_\phi \} \Big \} e^2.\\
\end{split}
\end{equation}
Observe, in the case of the sphere, $\abs{\nabla \rho}=1$, so if $E=S^2$, the formula immediately simplifies to

\be
\mathfrak{i}_{\mathrm{S^2}}^* \Big \{ -\triangle v \Big \} =   -\triangle_{\mathrm{S^2}} v  + \mathcal{E} \big ( v \big )
-  \mathfrak{i}_{\mathrm{S^2}}^* \Big \{ \mathcal{L}_{\mathrm{Y}} \big ( v \big ) \Big \} .
\ee

Next, from the definition of $\mathcal E$, we have 
\[
 \mathcal{E} = \mathfrak{i}_{\mathrm{E}}^* \Big \{ -\mathcal{L}_{\nabla \rho} \Big ( \frac{1}{|\nabla \rho|^2}  \mathcal{L}_{\nabla \rho}     \Big ) + \Big ( 1 - \frac{1}{|\nabla \rho|^2}     \Big )\mathcal{L}_{\nabla \rho} \Big \} +    \mathfrak{i}_{\mathrm{E}}^* \Big \{ \frac{1}{|\nabla \rho|^2}  \mathcal{L}_{\nabla \rho}  \Big \} .
\]
 We can further simplify to
\be
\mathfrak{i}_{\mathrm{S^2}}^* \Big \{ -\triangle v \Big \} =   -\triangle_{\mathrm{S^2}} v   +\mathfrak{i}_{\mathrm{S^2}}^* \Big \{ -\mathcal{L}_{\nabla \rho} \mathcal{L}_{\nabla \rho}  v  \Big\},
\ee
which to simplify further, we look at \eqref{Liedouble}.  For $S^2$, with $v=v^\phi(\phi,\theta) \partial_\phi+v^\theta(\phi,\theta)\partial_\theta$, we have
\begin{align}
 \mathfrak{i}_{\mathrm{S^2}}^* \Big \{   -\mathcal{L}_{\nabla \rho}  \mathcal{L}_{\nabla \rho}  v \Big\}
 =\mathfrak{i}_{\mathrm{S^2}}^* \Big \{ G_\phi \dd\phi + G_\theta \dd\theta \Big\},\label{Liedouble}
\end{align}
where 
\begin{align*}
G_\phi=-\partial_\rho \o_{\rho\phi},\quad
G_\theta=-\partial_\rho \o_{\rho\theta} .
\end{align*}
Now
\begin{align*}
G_\phi&=-\partial_\rho \o_{\rho\phi}=-\partial_\rho(\partial_\rho v_\phi-\partial_\phi v_\rho)= -\partial^2_\rho  (\rho^2v^\phi)=-2v^\phi=-2v_\phi\vert_{S^2} ,\\
G_\theta&=-\partial_\rho \o_{\rho\theta}=\partial_\rho(\partial_\rho v_\theta-\partial_\theta v_\rho)= -\partial^2_\rho  (\rho^2\sin^2\phi v^\theta)=-2\sin^2\phi v^\theta=-2v_\theta\vert_{S^2},
\end{align*}
so the result follows as needed.

\section{Existence of a smooth vector field tangential to the ellipsoid which is divergence free on $\R^3$ and on $E$}\label{divfree}

The defining function for $E$ which we are using is given by
\begin{equation}
\rho (x,y,z) = \Big ( \frac{x^2+ y^2 + a^2 z^2}{a^2} \Big )^{\frac{1}{2}},
\end{equation}
so then $E = \{ (x,y,z) \in \mathbb{R}^3 : \rho (x,y,z) = 1  \}$ .

Now, let $\o^{\sharp}$ be an arbitrary smooth vector field defined in the neighborhood $U$ of $E$, which is also tangential to $E$ and divergence free on $E$. 

Next, we add an extra component $v^{\rho} \partial_\rho$ to $\o^{\sharp}$ to create a smooth vector field $v^{\sharp}$
on $U\subset \mathbb{R}^3$.  The vector field $v^{\sharp}$ is defined by

\begin{equation}
v^{\sharp} = v^{\rho} \partial_\rho + \o^{\sharp} ,
\end{equation}
where $v^{\rho}$ is a smooth function on $U$, which we have to select so that $v^{\sharp}$
 is divergence free on $U\subset \R^3$ and on $E$.
Note that
\begin{equation}
\rho \partial_\rho = Y = x \partial_x + y \partial_y + z \partial_z ,
\end{equation}
thus:
\begin{equation}
\partial_\rho = \Big ( \frac{a^2}{x^2+ y^2 + a^2 z^2} \Big )^{\frac{1}{2}} \big (  x \partial_x + y \partial_y + z \partial_z \big ).
\end{equation}
This means that $\partial_\rho$ is a well-defined smooth vector field on $U$. Now, we insist that
$v^{\sharp}$ has to satisfy  
\begin{equation}\label{conditions}
\begin{split}
\dv_{\mathbb{R}^3} v^{\sharp} & = 0 ,\\
v^{\rho} \big |_E & = 0 .
\end{split}
\end{equation}
Observe
\begin{equation*}
\dv_{\mathbb{R}^3} \big ( v^{\rho} \partial_\rho \big ) = \partial_\rho v^{\rho} + \frac{2}{\rho} v^{\rho} .
\end{equation*}
Hence, we need
\begin{equation*}
0 = \dv_{\mathbb{R}^3} \big ( v^{\rho} \partial_\rho \big ) + \dv_{\mathbb{R}^3} (\o^\sharp) ,
\end{equation*}
which is equivalent to
\begin{equation*}
\partial_\rho \big ( \rho^2 v^\rho \big ) = - \rho^2 \dv_{\mathbb{R}^3} (\o^\sharp).
\end{equation*}

Now, we need to solve the above ordinary differential equation, with the initial condition
$v^{\rho} \big |_{E} = 0$.   For any point $(x_0 , y_0 , z_0) \in U$, the solution can be written as
\begin{equation}\label{vrhocomp}
v^{\rho} (x_0 , y_0 , z_0 ) = -\frac{1}{L_0^2} \int_1^{L_0} \tau^2 \dv_{\mathbb{R}^3} (\o^\sharp)   \dd \tau ,
\end{equation}
where 
\[
L_0 = \Big ( \frac{x_0^2 + y_0^2 + a^2 z_0^2 }{a^2} \Big )^{\frac{1}{2}} = \rho (x_0, y_0 , z_0 ) .
 \]

 By using the smooth function $v^{\rho}$ as defined by \eqref{vrhocomp}, it follows that $v^\sharp = v^{\rho} \partial_{\rho} + \o^\sharp$ must be a smooth function
on $\mathbb{R}^3 -\{O\}$ which is tangential to ellipsoid $E$, and which is divergence free both on $
\mathbb{R}^3$ and on the ellipsoid $E$. 

The last step is to let $f=\frac{v^\rho}{\rho-1}$ so that $v^\rho$ can be written as it was in Section \ref{proof}, i.e., as $f(\rho, \phi, \theta)(\rho-1)$.  It is a Calculus exercise to show that $f$ as defined is smooth.

The construction is thus completed.

\section{Asymptotic expansion of the main result }\label{power_exp}
In this section, we suppose $a>1$.  Then we consider the following parameter
\begin{equation}\label{eccent}
\mu = \Big ( \frac{a^2-1}{a^2} \Big )^{\frac{1}{2}} ,
\end{equation} 
which is the eccentricity of the cross section of the ellipsoid $E$ obtained by cutting through the ellipsoid by any plane containing the z-axis.  The eccentricity of a circle is $0$, so as $a\to 1^+$, we have $\mu\to 0^+$, in which case, the ellipsoid $E$ becomes a sphere of radius $1$. 

Then, we have the following useful information
\begin{equation}\label{extrauseful}
\begin{split}
\frac{1}{a^2} & = 1 - \mu^2 , \\
\frac{\lambda^2}{a^2} & = 1 - \mu^2 \sin^2 \phi ,\\
\partial_\phi \Big ( \frac{\lambda^2}{a^2}  \Big ) & = - 2 \mu^2 \sin \phi \cos \phi.  
\end{split}
\end{equation} 
Moreover, since $0< \mu < 1$, we have

\begin{equation}\label{geoSeries}
\frac{1}{1-\mu^2 \sin^2 \phi} = \sum_{k=0}^{\infty} \mu^{2k} \sin^{2k} \phi ,
\end{equation}  
and
\begin{equation}\label{diffgeoSeries}
\frac{1}{(1-\mu^2 \sin^2 \phi)^2} = \sum_{k=0}^{\infty} (k+1) \mu^{2k} \sin^{2k} \phi ,
\end{equation}
for all $\mu \in (0,1)$ and all $0 < \phi < \pi$.

In the following calculations, we let $\rho = 1$.  Again, for convenience, we recall the main formula
\[
\mathfrak{i}_{\mathrm{E}}^* \Big \{ -\triangle v \Big \} =   -\triangle_{\mathrm{E}}\Big ( \mathfrak{i}_{\mathrm{E}}^* v  \Big ) + \mathcal{E} \big ( v \big )
- \sqrt{K_{\mathrm{E}}}  \mathfrak{i}_{\mathrm{E}}^* \Big \{ \mathcal{L}_{\mathrm{Y}} \big ( v \big ) \Big \}
- 2 \sqrt[4]{K_{\mathrm{E}}} \Big (1- \frac{1}{|\nabla \rho |^2} \Big )\mathfrak{i}_{\mathrm{E}}^* \Big \{ (\mathcal L_{\nabla \rho} v)_\phi \} \Big \} e^2,\\
\]
and now focus on the terms that do not appear in the case of the sphere.

Recall that 
\[
\mathcal{E}(v) =  \mathcal{E}(v)_\phi \dd \phi + \mathcal{E}(v)_\theta \dd \theta , 
\]
where
\begin{equation}\label{repeateddata}
\begin{split}
\mathcal{E}(v)_\phi = & G_\phi - \partial_\phi \partial_\rho v^{\rho} + \Big ( \frac{\lambda^2 - a^2}{a^2} + \frac{1}{\lambda^2} \Big ) \o_{\rho \phi} ,\\
\mathcal{E}(v)_\theta = & G_\theta - \partial_\theta \partial_\rho v^\rho 
+ \Big ( 1 - \frac{a^2}{\lambda^2} + \frac{a^2}{\lambda^4} \Big ) 
\Big ( \frac{\lambda^2}{a^2} \o_{\rho \theta } + \o_{\phi \theta} \frac{1-a^2}{a^2} \sin \phi \cos \phi \Big ),
\end{split}
\end{equation}
where $G_\phi$ and $G_\theta$ are given in \eqref{G2}-\eqref{G1}.

Then, by direct calculations based on \eqref{extrauseful}- \eqref{diffgeoSeries}, we get
\begin{equation}\label{Gphi}
\begin{split}
G_{\phi} = & - \partial_\rho \o_{\rho \phi} 
+ \mu^2 \Big \{ \partial_\rho \o_{\rho \phi} \sin^2 \phi + 
\partial_\phi \big ( \sin \phi \cos \phi \o_{\rho \phi} \big ) + 2 \sin \phi \cos \phi \partial_\rho v^{\rho}  \Big \} +O(\mu^4)
\end{split}
\end{equation}
and 
\begin{equation}\label{Gtheta}
\begin{split}
G_\theta =  - \partial_\rho \o_{\rho \theta} + \mu^2 \Big \{ \sin^2 \phi \partial_\rho \o_{\rho \theta }  
+ \partial_\rho \Big ( \frac{\o_{\phi \theta}}{\rho } \Big ) \sin \phi \cos \phi 
+ \sin \phi \cos \phi \partial_\phi \o_{\rho \theta }  \Big \} +O(\mu^4)
\end{split}
\end{equation}

Now, we apply  \eqref{extrauseful}- \eqref{diffgeoSeries}, \eqref{Gphi} and \eqref{Gtheta} back to \eqref{repeateddata} to get

\begin{equation*}
\begin{split}
\mathcal{E}(v)_\phi = & -\partial_\rho \o_{\rho \phi } - \partial_\phi \partial_\rho v^{\rho} + \o_{\rho \phi} \\
& + \mu^2 \Big \{ \partial_\rho \o_{\rho \phi} \sin^2\phi  + \partial_\phi \big (\sin \phi \cos \phi \o_{\rho \phi} \big ) + 2 \sin \phi \cos \phi \partial_\rho v^{\rho} - \o_{\rho \phi} \Big \} +O(\mu^4),\\
\mathcal{E}(v)_\theta = &-\partial_\rho \o_{\rho \theta} - \partial_\theta \partial_\rho v^{\rho} + \o_{\rho \theta} \\
& + \mu^2 \Big \{ \sin^2 \phi \partial_\rho \o_{\rho \theta} + \partial_\rho \o_{\phi \theta} \sin \phi \cos \phi + \sin \phi \cos \phi \partial_\phi \o_{\rho \theta} - \o_{\rho \theta} -
 2 \o_{\phi \theta} \sin \phi \cos \phi \Big \} +O(\mu^4).
\end{split}
\end{equation*}

In exactly the same way, we can get

\begin{equation*}\label{veryextra1}
\begin{split}
-2 \sqrt[4]{K_{\mathrm{E}}} \Big ( 1 - \frac{1}{|\nabla \rho |^2}  \Big ) \mathfrak{i}_\mathrm{E}^* 
\Big \{ \Big ( \mathcal{L}_{\nabla \rho} v \Big )_\phi \Big \} e^2 =
2 \mu^2 \sin^2\phi \o_{\rho \phi} \dd \phi  +O(\mu^4), 
\end{split}
\end{equation*}

and 

\begin{equation*}
- \sqrt{K_{\mathrm{E}}} \mathfrak{i}_\mathrm{E}^* \Big \{ \mathcal{L}_Y v \Big \}
= \Big (-1 + \mu^2 \cos^2 \phi  \Big ) 
\Big ( \o_{\rho \phi} \dd \phi + \o_{\rho \theta} \dd \theta + \dd_{\mathrm{E}} v_\rho \Big )+O(\mu^4) .
\end{equation*}
We now combine the above results, and use the divergence free condition, to obtain 
 \begin{equation}\label{formula_asympt}
\begin{split}
& \mathfrak{i}_{\mathrm{E}}^* \Big \{ -\triangle v \Big \} 
=   -\triangle_{\mathrm{E}}\Big ( \mathfrak{i}_{\mathrm{E}}^* v  \Big ) -\partial_\rho \o_{\rho \phi } \dd\phi-\partial_\rho \o_{\rho \theta} \dd\theta\\
&\ +\mu^2 \Big \{ \sin^2\phi \partial_\rho \o_{\rho \phi}  + \cos^2\phi \o_{\rho\phi}+ \sin\phi\cos\phi\partial_\phi \omega_{\rho\phi}  \Big \} \dd\phi\\
&\ + \mu^2 \Big \{   \sin^2 \phi (\partial_\rho \o_{\rho \theta} - \o_{\rho \theta} )+  \sin \phi \cos \phi (\partial_\rho \o_{\phi \theta} +  \partial_\phi \o_{\rho \theta}   -
  2\o_{\phi \theta}  ) \Big \}\dd\theta\\
  &\  +O(\mu^4).
 \end{split}
\end{equation}
If we are interested in a differential form, then we can stop at the above formula. If we are interested in a formula involving explicitly the components of the vector field $v^\sharp$, then we can expand more as follows.  On $E$, from the divergence free condition, definition of the metric, and  $a^2-1\sim \mu^2$, we have

 \begin{equation}\label{formula_asympt2}
\begin{split}
& \mathfrak{i}_{\mathrm{E}}^* \Big \{ -\triangle v \Big \} 
=   -\triangle_{\mathrm{E}}\Big ( \mathfrak{i}_{\mathrm{E}}^* v  \Big ) -( 2v^\phi+4\partial_\rho v^\phi+\partial^2_\rho v^\phi)  \dd\phi- (2v^\theta+4\partial_\rho v^\theta+\partial_\rho^2 v^\theta) \dd\theta\\
&\quad-\mu^2\Big \{\partial^2_\rho v^\phi+ \partial_\rho v^\phi-\sin^2\phi (2v^\phi+4\partial_\rho v^\phi +2\partial^2_\rho v^\phi) \\
&\qquad\qquad - \sin\phi\cos\phi(2 \partial_\phi v^\phi+\partial_\rho \partial_\theta v^\theta+3\partial_\phi \partial_\rho v^\phi)  \Big \} \dd\phi\\
&\quad- \mu^2 \Big \{ \cos^2\phi (2v^\theta+4\partial_\rho v^\theta+\partial^2_\rho v^\theta) +(2\sin^2\phi\cos^2\phi-\sin^4\phi)(2v^\theta+\partial_\rho v^\theta) \\
&\qquad\qquad +2\sin\phi\cos\phi\left(\sin\phi \cos \phi \partial_\rho v^\theta+\sin^2\phi(\partial_\phi v^\theta+\partial_\rho\partial_\phi v^\theta)\right)\Big\} \dd\theta\\
&\quad +O(\mu^4).
  \end{split}
\end{equation}

\section*{Acknowledgements}
Chi Hin Chan is funded in part by a grant from the Ministry of Science and Technology of Taiwan (109-2115-M-009 -009 -MY2). This work was completed while Chi Hin Chan was working as a Center Scientist at the National Center for Theoretical Science of Taiwan R.O.C. 

Magdalena Czubak is funded in part by a grant from the Simons Foundation \# 585745.
 
Tsuyoshi Yoneda is funded in part by the JSPS Grants-in-Aid for Scientific
Research  17H02860, 18H01136, 18H01135 and 20H01819.

\bibliographystyle{plain}

\end{document}